\documentclass[11pt]{article}

\usepackage{amssymb}
\usepackage[mathscr]{euscript}
\usepackage{graphicx}

\newcommand{\dis}{\displaystyle}
\author{Fr\'ed\'eric H\'elein, ENS de Cachan}
\title{Isoperimetric Inequalities and Calibrations\footnote{published
in "Progress in Partial Differential Equations: the Metz surveys",
M. Chipot and I. Shafrir ed., Pitman Research Notes in Mathematics, Series 345, Longman (1996)}}
\begin{document}
\date{}
\newtheorem{defi}{Definition}
\newtheorem{theo}{Theorem}
\newtheorem{lemm}{Lemma}
\newtheorem{rema}{Remark}

\maketitle
\large
The subject of these Notes is a new proof, proposed in [4], of the classical isoperimetric inequality in the plane.
This proof is far from being the first one, and we may refer to [1] for a small review of some various proofs
which are known, and for further references. To my opinion, the interest of this proof is that it uses essentially integration
by parts and Stokes' formula in a simple manner, like in a calibration. Let us explain it briefly: consider a smooth domain
$\Omega$ of the plane $\Bbb{R}^2$, and let $x$, $y$ denote points of $\partial \Omega$, the boundary of $\Omega$. We denote by
$t_y$ the unit tangential vector to $\partial \Omega$ at $y$, such that, if $n_y$ is the exterior normal vector to $\partial \Omega$
at $y$, then $(n_y,t_y)$ forms a direct basis. We then let for $x \neq y$
$$V(y, t_y, x) = 2{\langle x-y, t_y\rangle \over |x-y|^2}(x-y) - t_y.$$
\begin{figure}[h]
\begin{center}
\includegraphics[scale=1]{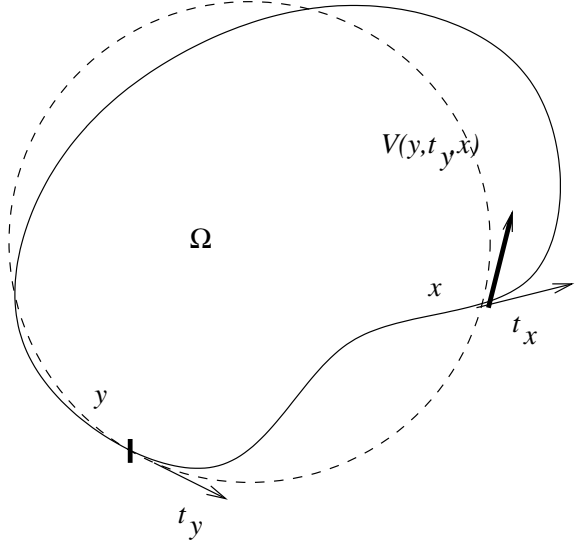}
\caption{construction of $V(y,t_y,x)$}
\end{center}
\end{figure}
We first fix $y$ and $t_y$ and we build the 1-form (see Figure 1)
$$\alpha = \langle V(y,t_y,x),dx\rangle .$$
We integrate $\alpha$ over $x\in \partial \Omega$. Using the fact that $V$ is of norm 1 everywhere, and Stokes' formula, we get
$$\begin{array}{lll}
|\partial \Omega |&\geq &\dis \int _{\partial \Omega} \alpha \\
&=&\dis \int _{x\in \Omega}d\alpha \\
&=&\dis \int_{x\in \Omega}2{det(y-x,t_y)\over |x-y|^2}dx^1\wedge dx^2.
\end{array}$$
Thus integrating this inequality over $y$, we get using Fubini's Theorem and Stokes's formula
$$\begin{array}{lll}
|\partial \Omega|^2&\geq&\dis \int _{y\in \partial \Omega} dl(y)\int _{x\in \Omega}2{det(y-x, t_y)\over |x-y|^2}dx^1\wedge dx^2\\
&=&\dis \int _{x\in \Omega}dx^1\wedge dx^2\int _{y\in \partial \Omega}2{det(y-x, d(y-x))\over |x-y|^2}\\
&=&\dis \int _{x\in \Omega}dx^1\wedge dx^2\int _{y\in \Omega}2d{det(y-x, d(y-x))\over |x-y|^2}\\
&=&\dis \int _{x\in \Omega}4\pi dx^1\wedge dx^2\\
&=&4\pi|\Omega|,
\end{array}$$
because $\dis d_y{det(y-x, d(y-x))\over |x-y|^2} = 4\pi \delta (y-x)dy^1\wedge dy^2$. An analogous demonstration gives also optimal
inequalities for domains in the sphere $S^2$ or the hyperbolic disc $H^2$, see [4]. 

This is a basic example of an important question in the calculus of variations: given a variational problem, and some particular
critical point, we would like to know whether this critical point is a global minimum of the functional involved. For many interesting
variational problems coming from Physics or Geometry, we know particular solutions which are candidates to be minimizers (harmonic maps,
Skyrme model, Ginzburg-Landau equation, Willmore problem...). But we have a very few tools for solving these questions. One of these tools
is the use of a null Lagrangian or a of calibration.

In the following we want to discuss in which sense our proof of the isoperimetric inequality works like a calibration. That is why we
first recall in Section 1 some classical notions about calibrations, and in particular we briefly expose the theory of null Lagrangians
in one variables. Although old (or maybe because of its age) this theory builded by Weierstrass, Mayer, Hilbert and later on generalized,
to the case of several variables, by Weyl, Caratheodory, Lepage and Boerner is not very popular. 

Then in the second Section we we come back to our proof of the isoperimetric inequality. It is interesting to observe how this proof has
some reminiscences from the theory of null Lagrangians, but also how it differs. This suggests that there is maybe a need for constructing
a more general theory for the null Lagrangians.

\section{Calibrations and null Lagrangians}

\subsection{Calibrations}
This is a kind of homological methods for proving that some "minimal" submanifolds, i.e. with vanishing mean curvature are area (or volume)
minimizing. The principle is the following. Let ${\cal N}$ be a Riemannian manifold and consider a $p$-dimensional submanifold ${\cal S}$ of
${\cal N}$. Assume that there exists a $p$-form $\alpha$ on ${\cal N}$ such that
\begin{equation}\label{5}\begin{array}{lllll}
(i)&|\alpha |&\leq&1&on\ {\cal N},\\
(ii)&|\alpha _{|{\cal S}}|&=&1& on\ {\cal S},\\
(iii)&d\alpha&=&0.&
\end{array}\end{equation}
Then $\alpha$ is called a $calibration$, and we will say that $\alpha$ $calibrates$ ${\cal S}$. As a consequence ${\cal S}$ is volume
minimizing, since for every submanifold ${\cal S}'$ which is homologous to ${\cal S}$ we have
$$|{\cal S}'| \geq \int _{{\cal S}'}\alpha = \int _{\cal S} \alpha = |{\cal S}|.$$
The advantage of this method is that it gives usually very short proofs. The defect is that it is difficult to find situations where
this principle works, and the instances which are known are generally founded by artisanal ways.

The first paper where this idea was applied is from Wirtinger [9] in 1936, but it seems that it was already discussed by
de Rham in unpublished lectures during the thirties. Federer used this method to prove that holomorphic submanifolds of a K\"ahler manifold
are minimizing using the K\"ahler form as a calibration. Surveys and many references may be found in [8] and [5].

\subsection{Null Lagrangians}
The only systematic and theoretical tool that we have at hand today is still the old theory that we present below. We consider here a very
simplified situation to explain the idea, a complete and detailled presentation may be founded in [6].

Let $(a, b)$ be some open interval of $\Bbb{R}$. We study variational problems for mappings from $(a, b)$ to $\Bbb{R}$. For this purpose
we choose a Lagrangian $L$, i.e. a function of the variables $(t, q, \dot q)$ in $(a, b)\times \Bbb{R} \times \Bbb{R}$. We assume that $L$
is sufficiently smooth ($C^2$ is enough) and that the Legendre condition
\begin{equation}\label{10}
{\partial ^2L\over \partial \dot q^2} \geq 0
\end {equation}
is true everywhere. We build on the space ${\cal E}=C^2((a, b),\Bbb{R})$ a functional ${\cal L}$ by
$${\cal L}(f) = \int_a^bL\left(t, f(t), {df\over dt}(t)\right)dt.$$
If $f_o\in {\cal E}$ is a solution of the Euler-Lagrange equation of ${\cal L}$:
$${d\over dt}\left[{\partial L\over \partial \dot q}\left(t, f_o(t), {df_o\over dt}(t)\right)\right]
= {\partial L\over \partial q}\left(t, f_o(t), {df_o\over dt}(t)\right),$$
it is natural to ask whether $f_o$ minimizes ${\cal L}$. This means, in order to make sense, that if we denote
$${\cal E}_o = \{f\in {\cal E}/f(a) = f_o(a), f(b) = f_o(b)\},$$
then for any $f$ in ${\cal E}_o$,
$${\cal L}(f) \geq {\cal L}(f_o).$$
In our situation an answer is
\begin{theo}. Let $f_o$ be a solution of the Euler-Lagrange equation. Assume that the two following conditions are true:
a) The Legendre condition (2) holds.
b) (Geometrical condition) There exists a foliation of $(a, b)\times \Bbb{R}$ by graphs of solutions of the Euler-Lagrange equation of
${\cal L}$, one of these solutions being $f_o$ itself.
Then $f_o$ is minimizing.
\end{theo}

Let us explain in more technical terms the second hypothesis. It means that there exists an interval $I$ of $\Bbb{R}$, and a map
$$\begin{array}{llll}
u:&I\times (a, b)&\longrightarrow &\Bbb{R}\\
&(s, t)&\longmapsto &u(s, t)
\end{array}$$
such that the map
$$\begin{array}{llll}
U:&I\times (a, b)&\longrightarrow &(a, b)\times \Bbb{R}\\
&(s, t)&\longmapsto &(t, u(s, t))
\end{array}$$
is a diffeomorphism. Moreover for some value of $s$, say $s=0$, $u(0,.)=f_o$, and for each $s$ fixed the function $t\longmapsto u(s,t)$ is a
solution of the Euler-Lagrange equation (see Figure 2).
\begin{figure}[h]
\begin{center}
\includegraphics[scale=1]{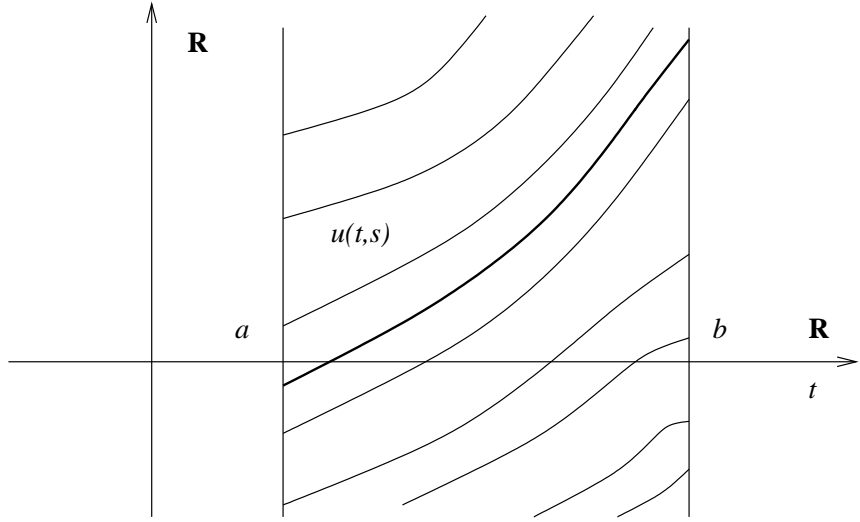}
\caption{foliation of $(a, b)\times \Bbb{R}$ by graphs of $t\longmapsto u(s, t)$}
\end{center}
\end{figure}

The proof of the Theorem is based on the construction of a null Lagrangian $\Lambda$ which "calibrates" $f_o$. 
\begin{defi} A null Lagrangian is a Lagrangian $\Lambda$ which satisfies one of the two following equivalent properties.
1) $\forall f \in {\cal E}$, $f$ is a solution of the Euler-Lagrange equation of $\Lambda$.
2) $\forall f\in {\cal E}, \int_a^b\Lambda \left(t, f(t), {df\over dt}(t)\right)dt$ depends uniquely on $f(a)$ and $f(b)$.
\end{defi}
Notice that it may be proved that any null Lagrangian $\Lambda$ can be constructed from a function $S:(a, b)\times \Bbb{R} \longrightarrow \Bbb{R}$ by
$$\Lambda (t, q, \dot q) = {\partial S\over \partial q}(t, q)\dot q +{\partial S\over \partial t}(t, q).$$

\noindent \underline{Proof.} We will build a Lagrangian $\Lambda$ such that
\begin{equation}\label{15}\begin{array}{lll}
(i)&\Lambda (t, q, \dot q) \leq L(t, q, \dot q)&\forall (t, q, \dot q)\in (a,b)\times \Bbb{R} \times \Bbb{R} ,\\
(ii)&\Lambda \left( t, f_o(t), {df_o\over dt}(t)\right) = L\left( t, f_o(t), {df_o\over dt}(t)\right) &\forall t\in (a,b),\\
(iii)&\Lambda \ is\ a\ null\ Lagrangian.&
\end{array}\end{equation}
We stress out the fact that the conditions (3) are the exact replicas of the conditions (1) for a calibration.

\noindent \underline{Step 1} For any $(t,q)$ in $(a,b)\times \Bbb{R}$, there exists a unique $s\in I$ such that the curve $\{(t, u(s, t))\}$
meets $(t,q)$. In other words
$$\exists ! s\in I, u(s, t) = q.$$
We then pose
$$\psi (t,q) = {\partial u\over \partial t}(s, t).$$
$\psi :(a,b)\times \Bbb{R} \longrightarrow \Bbb{R}$ is called a Mayer field.

\noindent \underline{Step 2} We build
$$\Lambda (t,q,\dot q) = \hat p(t,q,\psi (t,q))\dot q - E(t,q,\psi (t,q)),$$
where $\hat p:(a,b)\times \Bbb{R} \times \Bbb{R} \longrightarrow \Bbb{R}$ is given by
$$\hat p(t,q,\dot q) = {\partial L\over \partial \dot q}(t,q,\dot q)$$
(impulsion), and
$$E(t,q,\dot q) = {\partial L\over \partial \dot q}(t,q,\dot q)\dot q - L(t,q,\dot q)$$
is called the energy or Hamiltonian.

\noindent \underline{Step 3} We write down a Taylor formula
$$\begin{array}{lll}
L(t,q,\dot q) &=&\dis  L(t,q,\psi) + {\partial L\over \partial \dot q}(t,q,\psi)(\dot q -\psi)
+{\partial ^2L\over \partial \dot q ^2}(t,q,\psi +\theta (\dot q -\psi))(\dot q -\psi)^2\\
&=&\dis \Lambda (t,q,\dot q) +{\partial ^2L\over \partial \dot q ^2}(t,q,\psi +\theta (\dot q -\psi))(\dot q -\psi)^2,
\end{array}$$
where $\theta \in (0,1)$. Hence the Legendre condition (2) implies that
\begin{equation}\label{20}
L (t,q,\dot q) \geq \Lambda(t,q,\dot q).
\end{equation}

\noindent \underline{Step 4} We prove that $\Lambda$ is a null Lagrangian. For this purpose we will introduce Hamiltonian notations.
We consider the space $V =(a,b)\times \Bbb{R}^3$ equipped with coordinates $(t,q,e,p)$ and with the symplectic form
$\omega = dp\wedge dq - de\wedge dt$. We consider the Legendre transformation which exchange variables $(t,q,\dot q)$ and $(t,q,p)$
through the relation $p={\partial L\over \partial \dot q}(t,q,\dot q)$. For this we need that for any $(t,q,p)$ there exists a unique
$x$ solution to the equation
$$\hat p(t,q,x) = {\partial L\over \partial \dot q}(t,q,x) = p.$$
This is indeed the case if the Legendre condition is true. We then denote by $\hat q :(a,b)\times \Bbb{R}^2\longrightarrow \Bbb{R}$
the map defined by
\begin{equation}
{\partial L\over \partial \dot q}(t,q,\hat q(t,q,p)) = p.
\end{equation}
Lastly we define the Hamiltonian function
$$H(t,q,p) = p\hat q(t,q,p) - L(t,q,\hat q(t,q,p)).$$
And it is easy to prove using (5) that
\begin{equation}\left\{ \begin{array}{lll}
\dis {\partial H\over \partial q}(t,q,p) &=&\dis  -{\partial L\over \partial q}(t,q,\hat q(t,q,p))\\
\dis {\partial H\over \partial p}(t,q,p) &=& \hat q(t,q,p)
\end{array}\right.\end{equation}
Now let $f$ be a solution of the Euler-Lagrange equations for $L$, and let us pose $g(t) = \hat p(t,f(t),{df\over dt}(t))$. A simple computation
using (6) shows that the Hamilton equations
\begin{equation}\left\{ \begin{array}{lll}
\dis {df\over dt} &=& \dis {\partial H\over \partial p}(t,f(t),g(t))\\
\dis {dg\over dt} &=&\dis  -{\partial H\over \partial q}(t,f(t),g(t))
\end{array}\right.\end{equation}
hold.

We will apply these equations to the solutions involved in the foliation of $(a,b)\times \Bbb{R}$ related to the Mayer field. We define
$\Phi : I\times (a,b)\longrightarrow V$ by
$$\Phi (s,t) = (t,u(s,t),w(s,t),v(s,t)),$$
where
$$\begin{array}{lll}
v(s,t)&=&\hat p\left( t,u(s,t),{\partial u\over \partial t}(s,t)\right) \\
w(s,t)&=&H(t,u(s,t),v(s,t)).
\end{array}$$
We now prove that the pull-back image of $\omega$ by $\Phi$ vanishes. First we compute
$$\Phi ^*\omega = \left( {\partial v\over \partial t} {\partial u\over \partial s} -
 {\partial v\over \partial s} {\partial u\over \partial t} +  {\partial w\over \partial s}\right) dt\wedge ds.$$
But using the fact that for each $s$ fixed, the map $t\longmapsto (u(s,t),v(s,t))$ satisfies Hamilton's equations (7) we obtain that
$$ {\partial w\over \partial s} =  {\partial H\over \partial q} {\partial u\over \partial s}
+ {\partial H\over \partial p} {\partial v\over \partial s} = 
- {\partial v\over \partial t} {\partial u\over \partial s} +
 {\partial u\over \partial t} {\partial v\over \partial s},$$
which implies that $\Phi ^*\omega = 0$. Since $\Phi$ is an embedding this proves that the submanifold
${\cal S} := \Phi \left( I\times (a,b)\right)$ is a Lagrangian submanifold of $V$, i.e. that $\omega$ vanishes on ${\cal S}$.

By restricting to ${\cal S}$ the projection map $(t,q,e,p)\longmapsto (t,q)$, we see that ${\cal S}$ is diffeomorphic to $(a,b)\times \Bbb{R}$.
The inverse map is given by the "lifting" map 
$$\begin{array}{llll}
R:&(a.b)\times \Bbb{R} &\longrightarrow&{\cal S}\\
&(t,q)&\longmapsto&\left( t,q,H[t,q,p(t,q,\psi (t,q))] ,\hat  p(t,q,\psi (t,q))\right) .
\end{array}$$

Now we have a geometrical interpretation of the construction of $\Lambda$. We let $\alpha :=pdq -edt$, and we remark that
$$R^*\alpha = \hat p(t,q,\psi(t,q))dq - H\left( t,q,\hat p(t,q,\psi (t,q))\right) dt.$$
If we compare this last expression with the definition of $\Lambda$, we conclude that for any $f\in {\cal E}_o$, denoting
$\Gamma _f=\{(t,f(t))/t\in (a,b)\}$  the graph of $f$, we have
\begin{equation}\int _a^b\Lambda \left(t, f(t), {df\over dt}(t)\right)dt = \int _{\Gamma _f} R^*\alpha = \int _{R(\Gamma _f)} \alpha.\end{equation}
But we notice that
$$d\left( \alpha _{|{\cal S}}\right) = (d\alpha )_{|{\cal S}} = \omega _{|{\cal S}} = 0.$$
Thus using Stokes' formula we get
\begin{equation}\int _{R(\Gamma _f)} \alpha = \int _{R(\Gamma _{f_o})} \alpha.\end{equation}
And (8) and (9) imply that
\begin{equation}\int _a^b\Lambda \left(t, f(t), {df\over dt}(t)\right)dt = \int _a^b\Lambda \left(t, f_o(t), {df_o\over dt}(t)\right)dt.\end{equation}
\underline {Step 5} We conclude. For any $f\in {\cal E}_o$,
$${\cal L}(f) \geq \int _a^b\Lambda \left(t, f(t), {df\over dt}(t)\right)dt = \int _a^b\Lambda \left(t, f_o(t), {df_o\over dt}(t)\right)dt = {\cal L}(f_o).$$

Note that the above result can be generalized to variational problems with maps from a line segment $(a,b)$ to a vector space
$E$. But in this case the receipt does not work in general situations, and we need to require some further integrability conditions
on the foliation of $(a,b)\times E$ by graphs of solutions to the Euler-Lagrange equation. These conditions can be expressed geometrically
by assuming that the "lifting" of $(a,b)\times E$ into the symplectic space is a Lagrangian submanifold.

\section{The isoperimetric inequality}
By revisiting the method exposed in the introduction, we remark that $V$ plays the role of a Mayer field. A geometrical characterisation of
$V$ is the following: for each $(y,t_y,x) \in \Bbb{R} ^2\times S^1\times \Bbb{R}^2$ such that $x\neq y$, there exists a unique oriented
circle $C(y,t_y,x)$ in $\Bbb{R} ^2$ containing $x$ and $y$, and such that $t_y$ is tangent to this circle at $y$, and has a positive orientation.
Now $V(y,t_y,x)$ is nothing but the unit tangent vector to $C(y,t_y,x)$ at $x$, with a positive orientation - note that this geometrical construction
of $V(y,t_y,x)$ makes sense also on $S^2$ or $H^2$, that is how we can also prove optimal isoperimetric inequalities on $S^2$ and $H^2$, see [4]-.
Hence the Mayer field $V$ derives from a three parameters family of solutions: the set of all circles af the plane. But a major difference with the
classical theory is that the set of circles does not constitute a foliation of the plane. Instead it gives a foliation of
$\{(y,t,x)\in \Bbb{R} ^2\times S^1\times \Bbb{R} ^2\}$. Another difference also is that $V$ does not depends uniquely on a point $x$ of the
curve as expected in a classical theory but on two points plus a tangent vector.

We now propose tentative drafts for hypothetical theories which generalize the theory explained in Section 1, in order to reinterpret what
is happening in the isoperimetric inequality. We will adopt a more general framework: we define the ambiant space $X$ of the problem as follows.
If we are concerned with variational problems like in Section 1 (parametric problems) on maps from an open set $\Omega$ , into a vector space or
a Riemannian manifold $E$, $X$ is $\Omega \times E$, and to each map $f:\Omega \longrightarrow E$, we associate its graph in $X$.
In the case of a nonparametric problem, like minimal submanifolds, of a Riemannian manifold ${\cal N}$, then $X$ is ${\cal N}$.

\subsection{Use of families of foliations}
We would like here to construct null Lagrangians from more general objects than foliations of $X$. A simple generalisation is to start
with a family of foliations of $X$. For instance we introduce an auxiliary manifold $Y$, equipped with a measure $d\mu$ and we assume
that for any $y$ in $Y$, there exists a foliation of $X$ whose leaves are solutions to the variational problem. From this foliation we deduce a Mayer field and a Lagrangian using the classical receipt that we may represent by a differential form $\alpha _y$ in $X$. We then set 
$$\alpha = \int _Y \alpha _y d\mu (y).$$
And we could expect situations where the forms $\alpha _y$ are not closed, but the sum of all of them, $\alpha$ is closed. This is not enough
to have a calibration since we need also to find a submanifold ${\cal S}$ in $X$, such that after a renormalisation if necessary,
$$\begin{array}{llll}\
|\alpha _{y|{\cal S}}|&=&1&on\ {\cal S},\\
|\alpha |&\leq&1&elsewhere.
\end{array}$$
A similar mechanism appears in the isoperimetric inequality. Indeed whenever we fix $y \in \partial \Omega$, and $t_y$ tangent to $\partial \Omega$
at $y$, then the family $\{C(y,t_y,x)/x \in \Bbb{R}^2-\{y\}\}$ constitutes a foliation of $\Bbb{R}^2 -\{y\}$ by circles (see Figure 3). 
\begin{figure}[h]
\begin{center}
\includegraphics[scale=1]{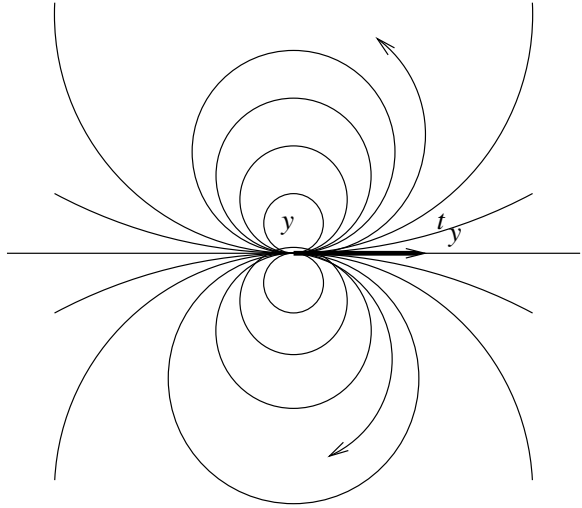}
\caption{foliation of $\Bbb{R}^2 -\{y\}$ by circles}
\end{center}
\end{figure}
And the one-form $\alpha _y = \langle V(y,t_y,x), dx\rangle $ derives from this foliation. Here $\alpha _y$ is not closed, and we do not want  $\alpha _y$ to be
closed because we want to estimate the length of $\partial \Omega$ in terms of the area of $\Omega$. But we would expect naively that
$\dis d\alpha _y = 4\pi {|\Omega |\over |\partial \Omega |}dx^1\wedge dx^2$, in such a way that
$$|\partial \Omega |^2\geq |\partial \Omega |\int_{\partial \Omega}\alpha _y = |\partial \Omega |\int _{\Omega}d\alpha _y = 4\pi |\Omega |.$$
This is of course not true. Instead when we sum all the $\alpha _y$'s
$$\alpha = {1\over |\partial \Omega |}\int _{\partial \Omega}\alpha _ydl(y),$$
we obtain this time that $|\alpha |\leq 1$ and
$$d\alpha = 4\pi {|\Omega |\over |\partial \Omega |}dx^1\wedge dx^2.$$
and we are lucky.

Of course it would be interesting to have other examples of calibrations constructed by this principle, particularly for variational problems
with several variables. This is because, according to my experience, in many interesting situations constructions by the classical theory with
several variables (Caratheodory, Weyl...) does not give efficient null Lagrangians. Notice also that when dealing with several variables the
problem is considerably more complicated than in one variable, since we have a too large choice of different possibilities and no way to choose
the best one, and also because of the heavy computations.

\subsection{Tensorial null Lagrangians}
The idea is to replace the ambiant space $X$ by $X^N$, for some $N\in \Bbb{N}$. For example, instead of looking for a null Lagrangian which works for
$${\cal L}(f) = \int _{\Omega}L\left(x,f(x),df(x)\right)dx.$$
where $f:\Omega \longrightarrow E$, do it for
$${\cal L}^N(f) = \int _{\Omega}...\int _{\Omega}L\left(x_1,f(x_1),df(x_1)\right)...L\left(x_N,f(x_N),df(x_N)\right)dx_1...dx_N.$$
Let assume $N=2$ for simplicity. We look for a lagrangian $\Lambda$ for maps from $\Omega \times \Omega$ into $E\times E$ such that first for all
$x,y \in \Omega$ and for any map $f:\Omega \longrightarrow E$,
\begin{equation}\begin{array}{l}
\Lambda \left[ (x,y),(f(x),f(y)),(df(x),df(y))\right]\\
\leq L(x,f(x),df(x)).L(y,f(y),df(y)),\end{array}\end{equation}
second, for a given map $f_o:\Omega \longrightarrow E$, equality holds in (11), and lastly $\Lambda$ is a null Lagrangian. Then it follows
that $f_o$ is minimizing.

In the framework of a minimal submanifolds it would give rise to the following. Let $\alpha$  be a $2p$-form on ${\cal M}\times {\cal M}$.
Assume that there exists a $p$-dimensional submanifold ${\cal S}$ of ${\cal N}$ such that
$$\begin{array}{lllll}
(i)&|\alpha |&\leq&1&on\ {\cal N}\times {\cal N}.\\
(ii)&|\alpha _{|{\cal S}\times {\cal S}}|&=&1& on\ {\cal S}\times {\cal S},\\
(iii)&d\alpha &=&0&.
\end{array}$$
Then ${\cal S}$ is volume minimizing.

Here also we do not have instances of this method but the isoperimetric inequality offers us a variant of this principle. This is made clearer
by the following presentation of our proof, which gives a symmetric role to the points $x$ and $y$. Consider for $n\in \Bbb{N}$ the tensor product
$(\bigwedge ^1\Bbb{R}^n)\otimes (\bigwedge ^1\Bbb{R}^n)$, where $(\bigwedge ^1\Bbb{R}^n)$ is the set of differential 1-forms on $\Bbb{R}^n$.
A tensor (or biform ?) $\beta$ in  $(\bigwedge ^1\Bbb{R}^n)\otimes (\bigwedge ^1\Bbb{R}^n)$ may be considered as a tensor on
$\{(x,y)\in \Bbb{R} ^n\times \Bbb{R} ^n\}$ and can be written as 
$$\beta = \sum_{i,j=1}^n\beta _{ij}(x, y)dx^i\otimes dy^j.$$
We define the differential
$$\begin{array}{cccc}
d_1:& (\bigwedge ^1\Bbb{R}^n)\otimes (\bigwedge ^1\Bbb{R}^n)&\longrightarrow& (\bigwedge ^2\Bbb{R}^n)\otimes (\bigwedge ^1\Bbb{R}^n)\\
&\beta &\longmapsto &\dis \sum_{i,j,k=1}^n{\partial \beta _{ij}\over \partial x^k}(x, y)dx^k\wedge dx^i\otimes dy^j,
\end{array}$$
and similarly $d_2: (\bigwedge ^1\Bbb{R}^n)\otimes (\bigwedge ^1\Bbb{R}^n)\longrightarrow (\bigwedge ^1\Bbb{R}^n)\otimes (\bigwedge ^2\Bbb{R}^n)$.

Now we set $n=2$ and
\begin{equation}\begin{array}{lll}
\alpha &=&\dis 2{\left[ (x^1-y^1)dx^1+(x^2-y^2)dx^2\right] \otimes \left[ (x^1-y^1)dy^1+(x^2-y^2)dy^2\right] \over |x-y|^2}\\
&-&\dis \left( dx^1\otimes dy^1+dx^2\otimes dy^2\right).
\end{array}\end {equation}
We observe that
\begin{equation}\begin{array}{llll}
(i)&|\alpha |&= &1\ on\ \Bbb{R}^2\times \Bbb{R}^2,\\
(ii)&|\alpha _{|C\times C}|&=&1\ for\ any\ circle\ C\subset \Bbb{R}^2,\\
(iii)&d_1d_2\alpha&=& 4\pi \delta (x-y)dx^1\wedge dx^2\otimes dy^1\wedge dy^2,
\end{array}\end{equation}
where $\delta$ is the Dirac distribution at the origin of $\Bbb{R}^2$. And thus for any $\Omega \subset \Bbb{R}^2$,
$$|\partial \Omega |^2\geq \int _{\partial \Omega}\int _{\partial \Omega}\alpha = \int _{\Omega}\int _{\Omega}d_1d_2\alpha = 4\pi |\Omega |.$$

$Remark\ 1.$ It is interesting to notice that we can also write $\alpha$ as

$\alpha = \sum_{i,j=1}^2 \alpha _{ij}(x,y)dx^i\otimes dy^j$, where
$$\left( \alpha _{ij}\right) = {1\over |x-y|^2}\left(
\begin{array}{cc}
(x^1-y^1)^2 - (x^2-y^2)^2&2(x^1-y^1)(x^2-y^2)\\
2(x^1-y^1)(x^2-y^2)&(x^2-y^2)^2 - (x^1-y^1)^2
\end{array}\right) .$$
This kind of matrix has been used by A.L. Bertozzi and P. Constantin in [2] for giving a new proof of a result from J.-Y. Chemin in [3]
about the regularity of a vortex patch in a 2-dimensional incompressible perfect fluid evolving according to Euler's equations
(another alternative proof was also given by P. Serfati in [7]). A vortex patch is modelised as an open set of $\Bbb{R}^2$ on which the
vorticity of the fluid is uniformly equal to one. Outside this open subset the vorticity vanishes. Here this matrix appears in an tricky integral
formulation of the evolution of the boundary of this open subset. The property of this matrix is that its "double curl" (as in (13) (iii)) gives
a Dirac mass at the origin. So it works like the kernel of a second order elliptic operator. But the advantage of $\alpha$ on
${1\over 2\pi}log(1/|x-y|)$ for instance is that it is bounded, which allows sharp estimates.

$Remark\ 2.$ We can also define in $ (\bigwedge ^1\Bbb{R}^3)\otimes (\bigwedge ^1\Bbb{R}^3)$ a tensor  generalizing obviously the definition as
in (12), i.e. denoting $z=x-y$,
$$\begin{array}{lll}
\alpha &=&\dis 2{\left( z^1dx^1 + z^2dx^2+z^3dx^3\right) \otimes \left( z^1dy^1+z^2dy^2+z^3dy^3\right)\over |z|^2}\\
&&-\left( dx^1\otimes dy^1+dx^2\otimes dy^2+dx^3\otimes dy^3\right).
\end{array}$$
We check that $\alpha$ satisfies
\begin{equation}\begin{array}{llll}
(i)&|\alpha |&= &1\ on\ \Bbb{R}^3\times \Bbb{R}^3,\\
(ii)&|\alpha _{|C\times C}|&=&1\ for\ any\ circle\ C\subset \Bbb{R}^3,
\end{array}\end{equation}
and lastly

\noindent $(iii)\ d_1d_2\alpha =$
$$4\pi \delta _2(x-y)\ [ 
dx^2\wedge dx^3\otimes dy^2\wedge dy^3 + dx^3\wedge dx^1\otimes dy^3\wedge dy^1 + dx^1\wedge dx^2\otimes dy^1\wedge dy^2]$$
$$-\ 4\left[ {(x^1-y^1)dx^2\wedge dx^3 +(x^2-y^2)dx^3\wedge dx^1+(x^3-y^3)dx^1\wedge dx^2\over |x-y|^2}\right.$$
$$\left.\otimes \ {(y^1-x^1)dy^2\wedge dy^3 +(y^2-x^2)dy^3\wedge dy^1+(y^3-x^3)dy^1\wedge dy^2\over |y-x|^2}\right].$$
Here $\delta _2$ refers to the Dirac distribution.
When restricting $\alpha$ to $S^2\times S^2$, we obtain
$$d_1d_2\alpha _{|S^2\times S^2} = \left(4\pi \delta(x-y) - 1\right)d\sigma (x)\otimes d\sigma (y),$$
where $d\sigma$ refers to the standard volume form on $S^2$. If we integrate as before $\alpha$ on $\partial \Omega \times \partial \Omega$,
where $\Omega$ is an open subset of $S^2$, it leads to the optimal isoperimetric inequality,
$$\begin{array}{lll}
|\partial \Omega |^2&\geq&\dis \int _{\partial \Omega}\int _{\partial \Omega}\alpha\\
&=&\dis 4\pi \int _{\Omega} d\sigma(y) - \int _{\Omega}d\sigma (x). \int _{\Omega}d\sigma (y)\\
&=&(4\pi -|\Omega )|)|\Omega |.
\end{array}$$
The proof of the isoperimetric inequality on  $H^2$ (i.e. $|\partial \Omega |^2 \geq (4\pi +|\Omega )|)|\Omega |$) can be done by a similar
computation in $\Bbb{R}^3$ equipped with a flat Minkowski metric in which we embedd isometrically $H^2$.

Fr\'ed\'eric H\'elein, 

Centre de Math\'ematiques et de Leurs Applications, 

Ecole Normale Sup\'erieure de Cachan, 

61 avenue du Pr\'esident Wilson 94235 Cachan Cedex, France.

Electronic mail: Frederic.Helein@cmla.ens-cachan.fr

\end{document}